\documentclass{article}
\usepackage{amsfonts}
\usepackage{amsmath}
\usepackage{amssymb}
\usepackage{graphicx}

\newtheorem{theorem}{Theorem}

\newtheorem{definition}[theorem]{Definition}

\newtheorem{lemma}[theorem]{Lemma}

\input{tcilatex}
\begin{document}

\title{\textbf{A Necessary and Sufficient Condition on the Weyl Manifolds
Admitting a Semi Symmetric Non-Metric Connection to be S-Concircular }}
\author{Assist. Prof. Dr. F\"{u}sun NURCAN BA\c{S}TAN \\
The Department of Mathematics,\\
Faculty of Science and Letters, \\
Marmara University, TURKEY. \\
funal@marmara.edu.tr}
\maketitle

\begin{abstract}
The object of this paper is to obtain the concircular curvature tensor of
the semi symmetric non-metric connection on the Weyl manifold and to give a
necessary and sufficient condition for a semi symmetric non-metric
connection to be S-concircular.
\end{abstract}

\textbf{Keywords: }Weyl manifolds, semi symmetric non-metric connection,
concircular curvature tensor, S-concircular connection.

\begin{quote}
\textbf{Mathematical Subject Classification :53A40}
\end{quote}

\section{\textbf{Introduction }}

An n-dimensional manifold which has a symmetric connection $\nabla $ and a
conformal metric tensor $g$ satisfying the compatibility condition%
\begin{equation*}
\nabla g=2\left( T\otimes g\right)
\end{equation*}%
where $T$ is a 1-form is called a Weyl space which is denoted by $W_{n}(g,T)$%
, (see \cite{Norden 1}). In local coordinates, the compatibility condition
is given by%
\begin{equation}
\nabla _{k}g_{ij}-2g_{ij}T_{k}=0  \tag{1.1}
\end{equation}%
where $T_{k}$ is a comlementary covariant vector field. Such a Weyl manifold
will be denoted by $W_{n}(g_{ij},T_{k}).$ If $T_{k}=0$ or $T_{k}$ is
gradient, a Riemannian manifold is obtained.

In \cite{Norden 1}, under the transformation of the metric tensor $g_{ij}$
in the form of%
\begin{equation}
\overset{\backsim }{g_{ij}}=\lambda ^{2}g_{ij}  \tag{1.2}
\end{equation}%
$T_{k}$ changes by%
\begin{equation*}
\overset{\backsim }{T_{k}}=T_{k}+\partial _{k}\left( \ln \lambda \right) ,
\end{equation*}%
where $\lambda $ is a scalar function defined on $W_{n}$.

The coefficients $\Gamma _{jk}^{i}$ of the symmetric connection $\nabla $ on
the Weyl manifold are defined by%
\begin{equation}
\Gamma _{jk}^{i}=\QTATOPD\{ \} {i}{jk}-g^{im}\left(
g_{mj}T_{k}+g_{mk}T_{j}-g_{jk}T_{m}\right)  \tag{1.3}
\end{equation}%
where $\QTATOPD\{ \} {i}{jk}$'s are the coefficients of the Levi-Civita
connection, (see \cite{Norden 1}).

In \cite{Norden 1}, the coefficients $\Gamma _{jk}^{i}$ and the curvature
tensor $R_{ijk}^{h}$ of the symmetric connection $\nabla $ change by%
\begin{equation}
\Gamma _{jk}^{i\ast }=\Gamma _{jk}^{i}+\delta _{j}^{i}P_{k}+\delta
_{k}^{i}P_{j}-g_{jk}P^{i}  \tag{1.4}
\end{equation}%
and%
\begin{equation}
R_{ijk}^{h\ast }=R_{ijk}^{h}+2\delta _{i}^{h}\nabla _{\left[ j\right.
}P_{\left. k\right] }+\delta _{k}^{h}P_{ij}-\delta
_{j}^{h}P_{ik}+g_{ij}g^{hr}P_{rk}-g_{ik}g^{hr}P_{rj}  \tag{1.5}
\end{equation}%
where $T_{i}-T_{i}^{\ast }=P_{i}$ and $P_{ij}=\nabla _{j}P_{i}-P_{i}P_{j}+%
\frac{1}{2}g_{ij}g^{kr}P_{k}P_{r}$ under conformal mapping $g_{ij}^{\ast
}=g_{ij}.$

The conformal curvature tensor $C_{ijk}^{h}$ and the concircular curvature
tensor $Z_{ijk}^{h}$ of the symmetric connection $\nabla $ on the Weyl
manifold are given by%
\begin{align}
C_{mijk}& =R_{mijk}-\frac{1}{n}g_{mi}R_{rjk}^{r}+\frac{1}{n-2}\left(
g_{mj}R_{ik}-g_{mk}R_{ij}-g_{ij}R_{mk}+g_{ik}R_{mj}\right)  \notag \\
& -\frac{1}{n\left( n-2\right) }\left(
g_{mj}R_{rki}^{r}-g_{mk}R_{rji}^{r}-g_{ij}R_{rkm}^{r}+g_{ik}R_{rjm}^{r}%
\right)  \tag{1.6} \\
& -\frac{R}{\left( n-1\right) \left( n-2\right) }\left(
g_{mj}g_{ik}-g_{mk}g_{ij}\right)  \notag
\end{align}%
and%
\begin{equation}
Z_{mijk}=R_{mijk}-\frac{R}{n\left( n-1\right) }(g_{mk}g_{ij}-g_{mj}g_{ik}), 
\tag{1.7}
\end{equation}%
where $R_{ijk}^{h}$, $R_{ij}$ and $R$ denote the curvature tensor, Ricci
tensor and scalar curvature tensor of $\nabla $, respectively, (see \cite%
{Miron 2}, \cite{Ozdeger 3}).

In \cite{Murgescu 4}, V.Murgescu defined the coefficients $\overline{\Gamma }%
_{jk}^{i}$ of a generalized connection $\overline{\nabla }$ on the Weyl
manifold by%
\begin{equation}
\overline{\Gamma }_{jk}^{i}=\Gamma _{jk}^{i}+a_{jkh}g^{hi}  \tag{1.8}
\end{equation}%
where%
\begin{equation}
a_{jkh}=g_{jr}\Omega _{kh}^{r}+g_{rk}\Omega _{jh}^{r}+g_{rh}\Omega _{jk}^{r}
\tag{1.9}
\end{equation}%
and $\Gamma _{jk}^{i}$'s are the coefficients of the symmetric connection $%
\nabla .$

By choosing%
\begin{equation*}
\Omega _{jk}^{i}=\delta _{j}^{i}a_{k}-\delta _{k}^{i}a_{j}
\end{equation*}%
in (1.9), the coefficients $\overline{\Gamma }_{jk}^{i}$'s of a semi
symmetric non-metric connection $\overline{\nabla }$ on the Weyl manifold
are obtained by%
\begin{equation}
\overline{\Gamma }_{jk}^{i}=\Gamma _{jk}^{i}+\delta _{k}^{i}S_{j}-g_{jk}S^{i}
\tag{1.10}
\end{equation}%
(see \cite{Unal F 5}). In (1.10), $S_{i}=-2a_{i}$ where $a_{i}$ is an
arbitrary covariant vector field.

The following results are also obtained in \cite{Unal F 5}:

The torsion tensor $T_{jk}^{i}$ with respect to the semi symmetric
connection $\overline{\nabla }$ is%
\begin{equation}
T_{jk}^{i}=\delta _{k}^{i}S_{j}-\delta _{j}^{i}S_{k}  \tag{1.11}
\end{equation}

The curvature tensor $\overline{R}(X,Y)Z$ of the semi symmetric non-metric
connection $\overline{\nabla }$ on the Weyl manifold is defined by%
\begin{equation*}
\overline{R}(X,Y)Z=\overline{\nabla }_{X}\overline{\nabla }_{Y}Z-\overline{%
\nabla }_{Y}\overline{\nabla }_{X}Z-\overline{\nabla }_{\left[ X,Y\right] }Z
\end{equation*}%
In local coordinates, above equation becomes%
\begin{equation}
\overline{R}_{ijk}^{h}=\partial _{j}\overline{\Gamma }_{ik}^{h}-\partial _{k}%
\overline{\Gamma }_{ij}^{h}+\overline{\Gamma }_{rj}^{h}\overline{\Gamma }%
_{ik}^{r}-\overline{\Gamma }_{rk}^{h}\overline{\Gamma }_{ij}^{r}  \tag{1.12}
\end{equation}

By means of (1.10) and (1.12), the relation between the curvature tensors $%
R_{ijk}^{h}$ and $\overline{R}_{ijk}^{h}$ of $\nabla $ and $\overline{\nabla 
}$, respectively, is obtained as%
\begin{equation}
\overline{R}_{ijk}^{h}=R_{ijk}^{h}+\delta _{k}^{h}S_{ij}-\delta
_{j}^{h}S_{ik}+g_{ij}g^{hr}S_{rk}-g_{ik}g^{hr}S_{rj}  \tag{1.13}
\end{equation}%
where%
\begin{equation}
S_{ij}=S_{i,j}-S_{i}S_{j}+\frac{1}{2}g_{ij}g^{kr}S_{k}S_{r}  \tag{1.14}
\end{equation}%
and $S_{i,j}$ denotes the covariant derivative of $S_{i}$ with respect to
the symmetric connection $\nabla $.

Transvecting (1.13) by $g_{mh}$ and contracting on the indices $h$ and $k$
in the same equation give 
\begin{equation}
\overline{R}%
_{mijk}=R_{mijk}+g_{mk}S_{ij}-g_{mj}S_{ik}+g_{ij}S_{mk}-g_{ik}S_{mj} 
\tag{1.15}
\end{equation}%
and%
\begin{equation}
\overline{R}_{ij}=R_{ij}+\left( n-2\right) S_{ij}+Sg_{ij}  \tag{1.16}
\end{equation}%
where$\ S=g^{mk}S_{mk}$, respectively.

The scalar curvatures $R$ and $\overline{R}$ of the connections $\nabla $
and $\overline{\nabla }$, respectively, are related by%
\begin{equation}
\overline{R}=R+2\left( n-1\right) S  \tag{1.17}
\end{equation}

The curvature tensor of the semi-symmetric connection $\overline{\nabla }$
has the following properties:

\begin{description}
\item[a] $\overline{R}_{mijk}+\overline{R}_{mikj}=0,$

\item[b] $\overline{R}_{mijk}+\overline{R}_{imjk}=4g_{mi}\nabla _{\lbrack
k}T_{j]},$

\item[c] $\overline{R}_{rjk}^{r}=R_{rjk}^{r}=2R_{\left[ kj\right] }=2n\nabla
_{\lbrack k}T_{j]},$

\item[d] $\overline{R}_{mijk}+\overline{R}_{mjki}+\overline{R}%
_{mkij}=2\left( g_{mi}\nabla _{\lbrack k}S_{j]}+g_{mj}\nabla _{\lbrack
i}S_{k]}+g_{mk}\nabla _{\lbrack j}S_{i]}\right) .$
\end{description}

The conformal curvature tensor $\overline{C}_{mijk}$ of $\overline{\nabla }$
is given by%
\begin{align}
\overline{C}_{mijk}& =\overline{R}_{mijk}-\frac{1}{n}g_{mi}\overline{R}%
_{rjk}^{r}+\frac{1}{n-2}\{g_{mj}\overline{R}_{ik}-g_{mk}\overline{R}%
_{ij}-g_{ij}\overline{R}_{mk}+g_{ik}\overline{R}_{mj}\}  \notag \\
& -\frac{1}{n(n-2)}\{g_{mj}\overline{R}_{rki}^{r}-g_{mk}\overline{R}%
_{rji}^{r}-g_{ij}\overline{R}_{rkm}^{r}+g_{ik}\overline{R}_{rjm}^{r}\} 
\tag{1.18} \\
& -\frac{\overline{R}}{\left( n-1\right) \left( n-2\right) }%
\{g_{mj}g_{ik}-g_{mk}g_{ij}\}.  \notag
\end{align}

The conformal curvature tensors $C_{mijk}$ and $\overline{C}_{mijk}$ of the
connections $\nabla $ and $\overline{\nabla }$ are related by

\begin{equation*}
\overline{C}_{mijk}=C_{mijk}
\end{equation*}

The projective curvature tensor $\overline{W}_{mijk}$ of $\overline{\nabla }$
is in the form of%
\begin{align}
\overline{W}_{mijk}& =\overline{R}_{mijk}+\frac{g_{mi}}{n+1}\left\{ \left( 
\overline{R}_{jk}-\overline{R}_{kj}\right) +2\left( n-1\right) \nabla
_{\lbrack j}S_{k]}\right\}  \tag{1.19} \\
& +\frac{1}{n^{2}-1}\left\{ g_{mj}\overline{H}_{ik}-g_{mk}\overline{H}%
_{ij}\right\}  \notag
\end{align}%
where%
\begin{equation*}
\overline{H}_{ij}=n\overline{R}_{ij}+\overline{R}_{ji}+2\left( n-1\right)
\nabla _{\lbrack j}S_{i]}\text{.}
\end{equation*}

The projective curvature tensors $W_{mijk}$ and $\overline{W}_{mijk}$ of the
connections $\nabla $ and $\overline{\nabla }$ are related by the equation 
\begin{equation*}
\overline{W}_{mijk}=W_{mijk}+\tfrac{2}{n+1}g_{mi}\nabla _{\lbrack j}S_{k]}+%
\tfrac{1}{n^{2}-1}\left( g_{mk}K_{ij}-g_{mj}K_{ik}\right)
+g_{ij}S_{mk}-g_{ik}S_{mj}
\end{equation*}%
where%
\begin{equation*}
K_{ij}=nS_{ij}+S_{ji}+\left( n+1\right) Sg_{ij}\text{.}
\end{equation*}

\section{\textbf{Weyl} \textbf{manifolds} \textbf{admitting a semi symmetric}%
\protect\linebreak \textbf{non-metric} \textbf{connection under concircular
mapping}}

Let $\ \sigma :(W_{n},g_{ij},T_{k},S_{k})\rightarrow (W_{n}^{\ast
},g_{ij}^{\ast },T_{k}^{\ast },S_{k}^{\ast })$ be a conformal mapping given
by $g_{ij}^{\ast }=g_{ij}$. In \cite{Unal F 5}, according to this mapping,
the coefficients $\overline{\Gamma }_{jk}^{i}$ and the curvature tensor $%
\overline{R}_{ijk}^{h}$ of the semi symmetric connection $\overline{\nabla }$
change by :%
\begin{equation}
\overline{\Gamma }_{jk}^{i^{\ast }}=\overline{\Gamma }_{jk}^{i}+\delta
_{j}^{i}P_{k}+\delta _{k}^{i}(P_{j}-Q_{j})-g_{jk}(P^{i}-Q^{i})  \tag{2.1}
\end{equation}%
where%
\begin{equation*}
P_{j}=T_{j}-T_{j}^{\ast },\text{ }Q_{j}=S_{j}-S_{j}^{\ast }
\end{equation*}%
and%
\begin{align}
\overline{R}_{ijk}^{h^{\ast }}& =\overline{R}_{ijk}^{h}+2\delta
_{i}^{h}\left( \nabla _{\lbrack j}P_{k]}+P_{[j}S_{k]}\right) +\delta
_{k}^{h}W_{ij}-\delta _{j}^{h}W_{ik}+g_{ij}g^{hr}W_{rk}  \notag \\
& -g_{ik}g^{hr}W_{rj}+2g^{sr}P_{s}Q_{r}\left( \delta _{j}^{h}g_{ik}-\delta
_{k}^{h}g_{ij}\right)  \tag{2.2}
\end{align}%
where%
\begin{equation*}
W_{ij}=\underline{P}_{ij}-\underline{Q}_{ij}+2P_{(i}Q_{j)},
\end{equation*}%
\begin{equation*}
\underline{P}_{ij}=P_{ij}-P_{i}S_{j},\text{ }P_{ij}=\nabla
_{j}P_{i}-P_{i}P_{j}+\frac{1}{2}g_{ij}g^{kr}P_{k}P_{r},
\end{equation*}%
\begin{equation*}
\underline{Q}_{ij}=Q_{ij}-Q_{i}S_{j},\text{ }Q_{ij}=\nabla
_{j}Q_{i}+Q_{i}Q_{j}-\frac{1}{2}g_{ij}g^{kr}Q_{k}Q_{r}.
\end{equation*}

Since a conformal mapping which $W_{ij}=\phi g_{ij}$ changes a geodesic
circle into a geodesic circle, it is called concircular mapping by means of 
\cite{Yano 6}.

Let $\sigma $ be a concircular mapping, that is, $\underline{P}_{ij}$ is
symmetric. Then, (2.2) can be rewritten as follows:%
\begin{equation}
\overline{R}_{ijk}^{h^{\ast }}=\overline{R}_{ijk}^{h}+2\left( \phi
-g^{sr}P_{s}Q_{r}\right) (\delta _{k}^{h}g_{ij}-\delta _{j}^{h}g_{ik}) 
\tag{2.3}
\end{equation}%
By contracting on h and k in (2.3),%
\begin{equation}
\overline{R}_{ij}^{\ast }=\overline{R}_{ij}+2\left( n-1\right) \left( \phi
-g^{sr}P_{s}Q_{r}\right) g_{ij}  \tag{2.4}
\end{equation}%
Transvecting (2.4) by $g^{ij^{\ast }}=g^{ij},$ yields%
\begin{equation}
\overline{R}^{\ast }=\overline{R}+2n\left( n-1\right) \left( \phi
-g^{sr}P_{s}Q_{r}\right)  \tag{2.5}
\end{equation}%
If the expression $2\left( \phi -g^{sr}P_{s}Q_{r}\right) =\frac{\overline{R}%
^{\ast }-\overline{R}}{n\left( n-1\right) }$ obtained from (2.5) is
substituted in (2.3),%
\begin{equation*}
\overline{R}_{ijk}^{h^{\ast }}=\overline{R}_{ijk}^{h}+\frac{\overline{R}%
^{\ast }-\overline{R}}{n\left( n-1\right) }(\delta _{k}^{h}g_{ij}-\delta
_{j}^{h}g_{ik})
\end{equation*}%
is arranged as%
\begin{equation*}
\overline{R}_{ijk}^{h^{\ast }}-\frac{\overline{R}^{\ast }}{n\left(
n-1\right) }(\delta _{k}^{h}g_{ij}^{\ast }-\delta _{j}^{h}g_{ik}^{\ast })=%
\overline{R}_{ijk}^{h}-\frac{\overline{R}}{n\left( n-1\right) }(\delta
_{k}^{h}g_{ij}-\delta _{j}^{h}g_{ik}).
\end{equation*}

If the \textit{concircular curvature tensor} $\overline{Z}_{ijk}^{h}$ is
defined by%
\begin{equation}
\overline{Z}_{ijk}^{h}=\overline{R}_{ijk}^{h}-\frac{\overline{R}}{n\left(
n-1\right) }(\delta _{k}^{h}g_{ij}-\delta _{j}^{h}g_{ik})\text{,}  \tag{2.6}
\end{equation}%
it is invariant under the concircular transformation, i.e. 
\begin{equation*}
\overline{Z}_{ijk}^{h^{\ast }}=\overline{Z}_{ijk}^{h}.
\end{equation*}

First transvecting (2.6) by $g_{mh}$ and then contracting on the indices $h$
and $k$ in (2.6), the equations 
\begin{equation}
\overline{Z}_{mijk}=\overline{R}_{mijk}-\frac{\overline{R}}{n\left(
n-1\right) }(g_{mk}g_{ij}-g_{mj}g_{ik})  \tag{2.7}
\end{equation}%
and%
\begin{equation}
\overline{Z}_{ij}=\overline{R}_{ij}-\frac{\overline{R}}{n}g_{ij}  \tag{2.8}
\end{equation}%
are obtained.

\begin{lemma}
\label{2.1} The concircular curvature tensor of the semi symmetric
connection $\overline{\nabla }$ has the following properties:
\end{lemma}

\begin{description}
\item[a] $\overline{Z}_{mijk}+\overline{Z}_{mikj}=0,$

\item[b] $\overline{Z}_{mijk}+\overline{Z}_{imjk}=4g_{mi}\nabla _{\lbrack
k}T_{j]},$

\item[c] $\overline{Z}_{rjk}^{r}=\overline{R}_{rjk}^{r},$

\item[d] $\overline{Z}_{mijk}+\overline{Z}_{mjki}+\overline{Z}_{mkij}=0.$
\end{description}

The concircular curvature tensors $Z_{ijk}^{h}$ and $\overline{Z}_{ijk}^{h}$
of $\nabla $ and $\overline{\nabla }$, respectively, are related by%
\begin{equation}
\overline{Z}_{ijk}^{h}=Z_{ijk}^{h}+\delta _{k}^{h}S_{ij}-\delta
_{j}^{h}S_{ik}+g_{ij}g^{mh}S_{mk}-g_{ik}g^{mh}S_{mj}-\frac{2}{n}S\left(
\delta _{k}^{h}g_{ij}-\delta _{j}^{h}g_{ik}\right)   \tag{2.9}
\end{equation}%
by substituting (1.7), (1.13) and (1.17) in (2.6).

Transvecting (2.9) by $g_{mh}$ and contracting on $h$ and $k$ in the same
equation give%
\begin{equation}
\overline{Z}%
_{mijk}=Z_{mijk}+g_{mk}S_{ij}-g_{mj}S_{ik}+g_{ij}S_{mk}-g_{ik}S_{mj}-\frac{2%
}{n}S\left( g_{mk}g_{ij}-g_{mj}g_{ik}\right)  \tag{2.10}
\end{equation}%
and%
\begin{equation}
\overline{Z}_{ij}=Z_{ij}+\left( n-2\right) S_{ij}-\frac{\left( n-2\right) }{n%
}g_{ij}S  \tag{2.11}
\end{equation}

\section{\textbf{Semi symmetric non-metric S-Concircular}\protect\linebreak 
\textbf{Connection}}

In \cite{Liang 7}, Liang defined semi symmetric recurrent metric connection
which is S-concircular on the Riemannian manifolds. In this paper, a semi
symmetric non-metric S-concircular connection on the Weyl manifold is
defined by as follows:

\begin{definition}
\label{4.1} If the semi symmetric non-metric connection $\overline{\nabla }$
satisfies the condition given by%
\begin{equation*}
S_{ij}=\nabla _{j}S_{i}-S_{i}S_{j}+\frac{1}{2}g_{ij}g^{rs}S_{r}S_{s}=\beta
g_{ij}
\end{equation*}%
where $\beta $ is a smooth function on the Weyl manifold, then it is called
S-concircular.
\end{definition}

Suppose that the concircular curvature tensors $Z_{mijk}$ and $\overline{Z}%
_{mijk}$ of symmetric and semi symmetric non-metric connections $\nabla $
and $\overline{\nabla }$, respectively, be the same. Then%
\begin{equation}
\overline{R}_{mijk}-\frac{\overline{R}}{n\left( n-1\right) }\left(
g_{mk}g_{ij}-g_{mj}g_{ik}\right) =R_{mijk}-\frac{R}{n\left( n-1\right) }%
\left( g_{mk}g_{ij}-g_{mj}g_{ik}\right) \text{.}  \tag{3.1}
\end{equation}%
By using (1.15) in (3.1), we get%
\begin{equation}
g_{mk}S_{ij}-g_{mj}S_{ik}+g_{ij}S_{mk}-g_{ik}S_{mj}=\frac{\overline{R}-R}{%
n(n-1)}(g_{mk}g_{ij}-g_{mj}g_{ik})  \tag{3.2}
\end{equation}%
By transvecting (3.2) by $g^{mk}$\thinspace , it is obtained as%
\begin{equation*}
(n-2)S_{ij}+Sg_{ij}=\frac{\overline{R}-R}{n}g_{ij}\text{ .}
\end{equation*}%
By (1.17),%
\begin{equation}
S_{ij}=\frac{\overline{R}-R}{2n(n-1)}g_{ij}  \tag{3.3}
\end{equation}%
which states that $\overline{\nabla }$ is S-concircular.

Conversely, suppose that $\overline{\nabla }$ is S-concircular. By using $%
S_{ij}=\beta g_{ij}$, from Definition \ref{4.1}, in (2.10), it is obtained as%
\begin{equation*}
\overline{Z}_{mijk}=Z_{mijk}.
\end{equation*}

In the view of the above results, we can state the following theorem:

\begin{theorem}
\label{4.2} The necessary and sufficient condition for the semi symmetric
non-metric connection $\overline{\nabla }$ to be S-concircular is that the
concircular curvature tensors $Z_{mijk}$ and $\overline{Z}_{mijk}$ of the
connections $\nabla $ and $\overline{\nabla }$, respectively, coincide.
\end{theorem}

\end{document}